\documentclass[9pt,article]{IEEEtran}
\usepackage{cite}
\usepackage{color}
\usepackage[pdftex]{graphicx}
\graphicspath{{fig/}{jpeg/}}
\usepackage[cmex10]{amsmath}
\usepackage{amssymb}
\usepackage{algorithm}
\usepackage{algorithmic}
\usepackage{amsmath,amssymb,lipsum}
\usepackage{hyperref}
\usepackage{comment}
\usepackage{graphicx}
\usepackage[font=small]{caption}
\usepackage{subcaption}
\usepackage{setspace}

\input{mysymbol.sty}
\usepackage{needspace}

\usepackage{theorem}
\newtheorem{thm}{Theorem}

\theoremstyle{definition}
\newtheorem{assumption}{Assumption}

\title{Network Newton}

\author{Aryan Mokhtari$^\dagger$, Qing Ling$^{\star}$ and Alejandro Ribeiro$^\dagger$ 
\\
{$^{\dagger}$Dept. of Electrical and Systems Engineering, University of Pennsylvania \\ 
    $^{\star}$Dept. of Automation, University of Science and Technology of China}
\thanks{Work in this paper is supported by ARO W911NF-10-1-0388, NSF CAREER CCF-0952867, and ONR N00014-12-1-0997. }
}

\begin{document}

\maketitle
\thispagestyle{empty}

\begin{abstract}



We consider minimization of a sum of convex objective functions where the components of the objective are available at different nodes of a network and nodes are allowed to only communicate with their neighbors. The use of distributed subgradient or gradient methods is widespread but they often suffer from slow convergence since they rely on first order information, which leads to a large number of local communications between nodes in the network. In this paper we propose the Network Newton (NN) method as a distributed algorithm that  incorporates second order information via distributed evaluation of approximations to Newton steps. We also introduce adaptive (A)NN in order to establish exact convergence. Numerical analyses show significant improvement in both convergence time and number of communications for NN relative to existing (first order) alternatives.

\end{abstract}


\section{Introduction}\label{sec_Introduction}

Distributed optimization algorithms are used to minimize a global cost function over a set of nodes in situations where the objective function is defined as a sum of a set of local functions. Consider a variable $\bbx\in\reals^p$ and a connected network containing $n$ agents each of which has access to a local function $f_i:\reals^p\to\reals$. The agents cooperate in minimizing the aggregate cost function $f:\reals^p\to\reals$ taking values $f(\bbx) := \sum_{i=1}^{n}f_i(\bbx)$. I.e., agents cooperate in solving the global optimization problem
\begin{equation}\label{original_optimization_problem1}
  \bbx^* \ :=\ \argmin_{\bbx}  f(\bbx) 
      \  =\ \argmin_{\bbx}\sum_{i=1}^{n} f_i(\bbx).
\end{equation}
Problems of this form arise often in, e.g., wireless systems \cite{Ribeiro10,Ribeiro12}, sensor networks \cite{Schizas2008-1,cRabbatNowak04}, and large scale machine learning \cite{Cevher2014}. There are different algorithms to solve \eqref{original_optimization_problem1} in a distributed manner. The most popular alternatives are decentralized gradient descent (DGD) \cite{Nedic2009,Jakovetic2014-1,YuanQing, Shi2014}, distributed implementations of the alternating direction method of multipliers \cite{Schizas2008-1,cQingRibeiroADMM14,BoydEtalADMM11,Shi2014-ADMM}, and decentralized dual averaging (DDA) \cite{Duchi2012}. A feature common to all of these algorithms is the slow convergence rate in ill-conditioned problems since they operate on first order information only. 

This paper considers Network Newton (NN), a method that relies on distributed approximations of Newton steps for the global cost function $f$ to accelerate convergence of the DGD algorithm. We begin this paper by introducing the idea that DGD solves a penalized version of \eqref{original_optimization_problem1} using gradient descent in lieu of solving the original optimization problem. To accelerate the convergence of gradient descent method for solving the penalty version of \eqref{original_optimization_problem1} we advocate the use of the NN algorithm. This algorithm relies on approximations to the Newton step of the penalized objective function by truncating the Taylor series of the exact Newton step (Section \ref{subsec_nn}). These approximations to the Newton step can be computed in a distributed manner with a level of locality controlled by the number $K$ of elements that are retained in the Taylor's series. When we retain $K$ elements in the series we say that we implement NN-$K$. We prove that for a fixed penalty coefficient lower and upper bounds on the Hessians of local objective functions $f_i$ are sufficient to guarantee at least linear convergence of NN-$K$ to the optimal arguments of penalized optimization problem (Theorem \ref{linear_convergence}). Further, We introduce an adaptive version of NN-$K$ (ANN-$K$) that uses an increasing penalty coefficient to achieve exact convergence to the optimal solution of \eqref{original_optimization_problem1} (Section \ref{sec:implement}). We study the advantages of NN-$K$ relative to DGD, both in terms of number of iterations and communications for convergence for solving a family of quadratic objective problems (Section \ref{sec:simulations}). 

\section{Problem formulation and Algorithm definition} \label{sec:problem}

The network that connects the agents is assumed symmetric and specified by the neighborhoods $\mathcal{N}_i$ that contain the list of nodes than can communicate with $i$ for $i=1,\ldots,n$. DGD is an established distributed method to solve \eqref{original_optimization_problem1} which relies on the introduction of local variables $\bbx_{i}\in\reals^p$ and nonnegative weights $w_{ij}\geq0$ that are not null if and only if  $j=i$ or if $j\in \mathcal{N}_i$. Letting $t\in\naturals$ be a discrete time index and $\alpha$ a given stepsize, DGD is defined by the recursion
\vspace{-1mm}
\begin{equation}\label{gd_iteration}
\vspace{-1mm}
   \bbx_{i,t+1} 
      = \sum_{j=1}^{n} w_{ij}\bbx_{j,t}-\alpha\nabla f_i({\bbx_{i,t}}), 
      \qquad i=1,\ldots,n.
\end{equation}
Since $w_{ij}=0$ when $j\neq i$ and $j\notin \mathcal{N}_i$, it follows from \eqref{gd_iteration} that each agent $i$ updates its estimate $\bbx_i$ of the optimal vector $\bbx^*$ by performing an average over the estimates $\bbx_{j,t}$ of its neighbors $j\in \mathcal{N}_i$ and its own estimate $\bbx_{i,t}$, and descending through the negative local gradient $-\nabla f_i(\bbx_{i,t})$. Note that weights $w_{ij}$ that nodes assign to each other form a weight matrix $\bbW\in\reals^{n\times n}$ that is symmetric and row stochastic. It is also customary to require the rank of $\bbI-\bbW$ to be $n-1$ so that $\text{null}(\bbI-\bbW)=\text{span}(\bbone)$. If the two assumptions $\bbW^T=\bbW$ and $\text{null}(\bbI-\bbW)=\bbone$ are true, it is possible to show that \eqref{gd_iteration} approaches the solution of \eqref{original_optimization_problem1} in the sense that $\bbx_{i,t}\approx\bbx^*$ for all $i$ and large $t$, \cite{Nedic2009}. 

To rewrite \eqref{gd_iteration} define the matrix $\bbZ:= \bbW \otimes \bbI\in\reals^{np\times np}$ as the Kronecker product of weight matrix $\bbW\in\reals^{n\times n}$ and the identity matrix $\bbI\in\reals^{p\times p}$. Further, we introduce vectors $\bby:= \left[ \begin{matrix} 
		\bbx_{1}; \dots ; \bbx_{n}
			\end{matrix} \right] \in\reals^{np}$ that concatenates the local vectors $\bbx_{i}$, and  vector $\bbh(\bby):= \left[ \begin{matrix} 
		\nabla f_{1}(\bbx_{1}); \dots ; \nabla f_{n}(\bbx_{n})
			\end{matrix} \right]\in\reals^{np}$ which concatenates the gradients of the local functions $f_i$ taken with respect to the local variable $\bbx_i$. It is then ready to see that \eqref{gd_iteration} is equivalent to
\begin{equation}\label{new_formulation2}
{   \bby_{t+1}\ =\ \bbZ \bby_t - \alpha\bbh(\bby_t)
             \ =\ \bby_{t} - \big[  (\bbI -\bbZ) \bby_t +\alpha \bbh(\bby_t)\big],}
\end{equation}
where in the second equality we added and subtracted $\bby_t$ and regrouped terms. Inspection of \eqref{new_formulation2} reveals that the DGD update formula at step $t$ is equivalent to a (regular) gradient descent algorithm being used to solve the program
\begin{equation}\label{centralized_opt_problem}
    \bby^* := \argmin\ F(\bby) 
           := \min \frac{1}{2}\ \bby^{T}(\bbI -\bbZ)\ \bby 
              + \alpha\sum_{i=1}^{n} f_i(\bbx_i).
\end{equation}
Observe that it is possible to write the gradient of $F(\bby)$ as
\begin{equation}\label{eqn_gradient_definition}
    \bbg_t \ :=\ \nabla F(\bby_t) 
           \  =\ (\bbI -\bbZ) \bby_t +\alpha \bbh(\bby_t),
\end{equation}
in order to write \eqref{new_formulation2} as $\bby_{t+1}=\bby_{t}- \bbg_t$ and conclude that DGD descends along the negative gradient of $F(\bby)$ with unit stepsize. The expression in \eqref{gd_iteration} is just a local implementation of \eqref{eqn_gradient_definition} where node $i$ implements the descent $\bbx_{i,t+1}=\bbx_{i,t}- \bbg_{i,t}$ where $\bbg_{i,t}$ is the $i$th element of the gradient $\bbg_{t}=[\bbg_{i,t};\ldots;\bbg_{i,t}]$. Node $i$ can compute the local gradient
\begin{equation}\label{local_gradient}
\bbg_{i,t}=(1-w_{ii})\bbx_{i,t} - \sum_{j\in \mathcal{N}_i} w_{ij} \bbx_{j,t}+\alpha \nabla f_{i}(\bbx_{i,t}).
\end{equation}

Notice that since we know that the null space of $\bbI-\bbW$ is $\text{null}(\bbI-\bbW)=\text{span}(\bbone)$ and that $\bbZ= \bbW \otimes \bbI$, we obtain that the span of $\bbI-\bbZ$ is $\text{null}(\bbI-\bbZ)=\text{span}(\bbone\otimes\bbI)$. Thus, we have that $(\bbI-\bbZ)\bby=\bb0$ holds if and only if $\bbx_1=\dots=\bbx_n$. Since the matrix $\bbI-\bbZ$ is positive semidefinite -- because it is stochastic and symmetric --, the same is true of the square root matrix $({\bbI-\bbZ})^{1/2}$. Therefore, we have that the optimization problem in \eqref{original_optimization_problem1} is equivalent to the optimization problem 
\begin{align}\label{original_optimization_new_notation}
   \tby^*: = \argmin_{\bbx}\ \sum_{i=1}^{n}\ f_i(\bbx_i), \qquad
   \text{s.t.} \quad ({\bbI-\bbZ})^{1/2}  \bby =\bb0.
\end{align} 
Indeed, for $\bby=[\bbx_1;\ldots;\bbx_n]$ to be feasible in \eqref{original_optimization_new_notation} we must have $\bbx_1=\dots=\bbx_n$ because $\text{null}[(\bbI-\bbZ)^{1/2}]=\text{span}(\bbone\otimes\bbI)$ as already argued. When restricted to this feasible set the objective $\sum_{i=1}^{n}\ f_i(\bbx_i)$ of \eqref{original_optimization_new_notation} is the same as the objective of \eqref{original_optimization_problem1} from where it follows that a solution $\tby^*=[\tbx_1^{*};\ldots;\tbx_n^{*}]$ of \eqref{original_optimization_new_notation} is such that $\tbx_i^{*}=\tbx^*=\bbx^*$ for all $i$, i.e. $\tby^*=[\bbx_1^{*};\ldots;\bbx_n^{*}]$. The unconstrained minimization in \eqref{centralized_opt_problem} is a penalty version of \eqref{original_optimization_new_notation}. The penalty function associated with the constraint $({\bbI-\bbZ})^{1/2}  \bby =\bb0$ is the squared norm $(1/2)\|({\bbI-\bbZ})^{1/2}  \bby \|^2$ and the corresponding penalty coefficient is $1/\alpha$. Inasmuch as the penalty coefficient $1/\alpha$ is sufficiently large, the optimal arguments $\bby^*$ and $\tby^*$ are not too far apart. In this paper we exploit the reinterpretation of \eqref{new_formulation2} as a method to minimize \eqref{centralized_opt_problem} to propose an approximate Newton algorithm that can be implemented in a distributed manner. We explain this algorithm in the following section.

%
\subsection{Network Newton}\label{subsec_nn}

Instead of solving \eqref{centralized_opt_problem} with a gradient descent algorithm as in DGD, we can solve \eqref{centralized_opt_problem} using Newton's method. To implement Newton's method we need to compute the Hessian $\bbH_t:=\nabla^2 F(\bby_t)$ of $F$ evaluated at $\bby_t$ so as to determine the Newton step $\bbd_t:=-\bbH_t^{-1}\bbg_t$. Start by differentiating twice in \eqref{centralized_opt_problem} in order to write $\bbH_t$ as
\begin{equation}\label{Hessian}
   \bbH_t := \nabla^2 F(\bby_t) = \bbI-\bbZ +\alpha \bbG_t,
\end{equation}
where the matrix $\bbG_t \in \reals^{np\times np}$ is a block diagonal matrix formed by blocks $\bbG_{ii,t}\in \reals^{p\times p}$ containing the Hessian of the $i$th local function,
\begin{equation}\label{G_form}
   \bbG_{ii,t} = \nabla^2 f_i(\bbx_{i,t})  .
\end{equation}
It follows from \eqref{Hessian} and \eqref{G_form} that the Hessian $\bbH_t$ is block sparse with blocks $\bbH_{ij,t}\in \reals^{p\times p}$ having the sparsity pattern of $\bbZ$, which is the sparsity pattern of the graph. {The diagonal blocks are of the form $\bbH_{ii,t}=(1-w_{ii})\bbI +  \alpha \nabla^2 f_i(\bbx_{i,t}) $ and the off diagonal blocks are not null only when $j\in \mathcal{N}_i$ in which case $\bbH_{ij,t}=w_{ij}\bbI$.} 

While the Hessian $\bbH_t$ is sparse, the inverse $\bbH_t$ is not. It is the latter that we need to compute the Newton step $\bbd_t:=-\bbH_t^{-1}\bbg_t$. To overcome this problem we split the diagonal and off-diagonal blocks of $\bbH_t$ and rely on a Taylor's expansion of the inverse. To be precise, write $\bbH_t=\bbD_t - \bbB$ where the matrix $\bbD_t$ is defined as 
\begin{equation}\label{diagonal_matrix}
   \bbD_t := \alpha \bbG_t + 2\ ( \bbI  -  \diag(\bbZ))
          := \alpha \bbG_t + 2\ ( \bbI  -  \bbZ_{d}),
\end{equation}
{in the second equality we defined $\bbZ_{d}:=\diag(\bbZ)$ for future reference.} Since the diagonal weights must be $w_{ii}<1$, the matrix $\bbI  -  \bbZ_{d}$ is positive definite. The same is true of the block diagonal matrix $\bbG_t$ because the local functions are assumed strongly convex. Therefore, the matrix $\bbD_t$ is block diagonal and positive definite. The $i$th diagonal block $\bbD_{ii,t}\in\reals^p$ of $\bbD_t$ can be computed and stored by node $i$ as $\bbD_{ii,t}= \alpha \nabla^2 f_{i}(\bbx_{i,t}) + 2(1-w_{ii})\bbI $. To have $\bbH_t=\bbD_t - \bbB$ we must define $\bbB:=\bbD_t-\bbH_t$.  Considering the definitions of $\bbH_t$ and $\bbD_t$ in \eqref{Hessian} and \eqref{diagonal_matrix}, it follows that 
\begin{equation}\label{non_diagona_matrix}
   \bbB =  \bbI - 2\bbZ_{d} +\bbZ.
\end{equation}
Observe that $\bbB$ is independent of time and depends on the weight matrix $\bbZ$ only. As in the case of the Hessian $\bbH_t$, the matrix $\bbB$ is block sparse with with blocks $\bbB_{ij}\in \reals^{p\times p}$ having the sparsity pattern of $\bbZ$, which is the sparsity pattern of the graph. Node $i$ can compute the diagonal blocks $\bbB_{ii}=(1-w_{ii})\bbI$ and the off diagonal blocks $\bbB_{ij}=w_{ij}\bbI$ using the local information about its own weights only.

Proceed now to factor $\bbD_t^{1/2}$ from both sides of the splitting relationship to write $\bbH_t = \bbD_t ^{{1}/{2}} (  \bbI - \bbD_t ^{{1}/{2}}\bbB\bbD_t ^{{1}/{2}} )^{-1} \bbD_t^{{1}/{2}}$. When we consider the Hessian inverse $\bbH^{-1}$, we can use the Taylor series  $(\bbI-\bbX)^{-1}= \sum_{j=0}^{\infty} \bbX^{j}$ with $\bbX=\bbD_t^{-{1}/{2}}  \bbB   \bbD_t^{-{1}/{2}}$ to write  
\begin{equation}\label{exact_Hessian_inverse}
   \bbH_t^{-1} = \bbD_t^{-1/2} 
                \sum_{k=0}^{\infty} \left(\bbD_t^{-1/2}  
                  \bbB   \bbD_t^{-1/2}\right)^{k}\ \bbD_t^{-1/2}.
\end{equation}
Observe that the sum in \eqref{exact_Hessian_inverse} converges if the absolute value of all the eigenvalues of matrix $\bbD^{-{1}/{2}}  \bbB   \bbD^{-{1}/{2}} $ are strictly less  than 1. This result is proven in \cite{NNICASSP}. 
%
\begin{algorithm}[t]
\caption{Network Newton-$K$ method at node $i$}\label{algo_NN1} 
\begin{algorithmic}[1] {
\REQUIRE  Initial iterate $\bbx_{i,0}$. 
\FOR {$t=0,1,2,\ldots$}
   \STATE Exchange iterates $\bbx_{i,t}$ with neighbors $j\in \mathcal{N}_i$.
   \STATE Gradient:    
          $\displaystyle{
          \bbg_{i,t} = (1-w_{ii})\bbx_{i,t} 
                       - \sum_{j\in \mathcal{N}_i} w_{ij} \bbx_{j,t}
                       +\alpha \nabla f_{i}(\bbx_{i,t}).}$  
   \STATE Compute NN-0 descent direction $\bbd_{i,t}^{(0)}=-\bbD_{ii,t}^{-1}\bbg_{i,t}$\\ 
   \FOR  {$k=  0, \ldots, K-1$ } 
      \STATE Exchange local elements $\bbd_{i,t}^{(k)}$ of the NN-$k$ step with neighbors
      \STATE NN-$(k+1)$ step:
             $\displaystyle{  
             \bbd_{i,t}^{(k+1)} = \bbD_{ii,t}^{-1}\!
             			\left[
                                   \sum_{j\in \mathcal{N}_i,j=i}\!\!\!\!\!\bbB_{ij} \bbd_{j,t}^{(k)} 
                                    - \bbg_{i,t}\right]}$.        
   \ENDFOR
          %
          \STATE Update local iterate: 
          $\displaystyle{\bbx_{i,t+1}=\bbx_{i,t} +\eps\ \bbd_{i,t}^{(K)}}$.
\ENDFOR}
\end{algorithmic}\end{algorithm}
%
Network Newton (NN) is defined as a family of algorithms that rely on truncations of the series in \eqref{exact_Hessian_inverse}. The {$K$th} member of this family, NN-$K$ considers the first $K+1$ terms of the series to define the approximate Hessian inverse
\begin{equation}\label{Hessian_inverse_approximation}
   \hbH_t^{(K)^{-1}} :=     \bbD_t^{-1/2}  \  \sum_{k=0}^{K} \left(  \bbD_t^{-1/2}  \bbB   \bbD_t^{-1/2}            \right)^{k}     \ \bbD_t^{-1/2}.
\end{equation}
NN-$K$ uses the approximate Hessian $\hbH_t^{(K)^{-1}}$ as a curvature  correction matrix that is used in lieu of the exact Hessian inverse $\bbH^{-1}$ to estimate the Newton step. I.e., instead of descending along the Newton step $\bbd_t:=-\bbH_t^{-1}\bbg_t$ we descend along the NN-$K$ step $\bbd_t^{(K)}:=-\hbH_t^{(K)^{-1}}\bbg_t$, which we intend as an approximation of $\bbd_t$. Using the explicit expression for $\hbH_t^{(K)^{-1}}$ in \eqref{Hessian_inverse_approximation} we write the NN-$K$ step as
\begin{equation}\label{Hessian_approximation_iteration}
\bbd_t^{(K)} = -\  \bbD_t^{-1/2}  \  \sum_{k=0}^{K} \left(  \bbD_t^{-1/2}  \bbB   \bbD_t^{-1/2}            \right)^{k}     \ \bbD_t^{-1/2}\ \bbg_t,
\end{equation}
where, we recall, the vector $\bbg_t$ is the gradient of objective function $F(\bby)$ defined in \eqref{eqn_gradient_definition}. The NN-$K$ update formula can then be written as
\begin{equation}\label{update_formula_NN}
   \bby_{t+1}=\bby_t+\eps\  \bbd_t^{(K)}.
\end{equation}
The algorithm defined by recursive application of \eqref{update_formula_NN} can be implemented in a distributed manner because the truncated series in \eqref{Hessian_inverse_approximation} has a local structure controlled by the parameter $K$. To explain this statement better define the components $\bbd^{(K)}_{i,t}\in\reals^p$ of the NN-$K$ step $\bbd^{(K)}_{t}=[\bbd^{(K)}_{1,t};\ldots;\bbd^{(K)}_{n,t}]$. A distributed implementation of \eqref{update_formula_NN} requires that node $i$ computes $\bbd^{(K)}_{i,t}$ so as to implement the local descent $\bbx_{i,t+1}=\bbx_{i,t} + \eps\bbd^{(K)}_{i,t}$. The step components $\bbd^{(K)}_{i,t}$ can be computed through local computations. To see that this is true first note that considering the definition of the NN-$K$ descent direction in \eqref{Hessian_approximation_iteration} the sequence of NN descent directions satisfies 
\begin{equation}
   \bbd_t^{(k+1)} = \bbD_t^{-1}\bbB \bbd_t^{(k)} -\bbD_t^{-1}\bbg_t
                  = \bbD_t^{-1}\left(\bbB \bbd_t^{(k)} - \bbg_t \right).
\end{equation}
Then observe that since the matrix $\hbB$ has the sparsity pattern of the graph, this recursion can be decomposed into local components
\begin{equation}\label{local_descent}
  \bbd_{i,t}^{(k+1)} 
      = \bbD_{ii,t}^{-1}\bigg(
              \sum_{j\in \mathcal{N}_i,j=i} \bbB_{ij} \bbd_{j,t}^{(k)} 
              - \bbg_{i,t}\bigg),
\end{equation}
The matrix $\bbD_{ii,t}=\alpha \nabla^2 f_{i}(\bbx_{i,t}) + 2(1-w_{ii})\bbI$ is stored and computed at node $i$. The gradient component $\bbg_{i,t}=(1-w_{ii})\bbx_{i,t} - \sum_{j\in \mathcal{N}_i} w_{ij} \bbx_{j,t}+\alpha \nabla f_{i}(\bbx_{i,t})$ is also stored and computed at $i$. Node $i$ can also evaluate the values of the matrix blocks $\bbB_{ij}=w_{ij}\bbI$. Thus, if the NN-$k$ step components $\bbd_{j,t}^{(k)}$ are available at neighboring nodes $j$, node $i$ can then determine the  NN-$(k+1)$ step component $\bbd_{i,t}^{(k+1)}$ upon being communicated that information.

The expression in \eqref{local_descent} represents an iterative computation embedded inside the NN-$K$  recursion in \eqref{update_formula_NN}. For each time index $t$, we compute the local component of the NN-$0$ step $\bbd_{i,t}^{(0)}=-\bbD_{ii,t}^{-1}\bbg_{i,t}$. Upon exchanging this information with neighbors we use \eqref{local_descent} to determine the NN-$1$ step components $\bbd_{i,t}^{(1)}$. These can be exchanged and plugged in \eqref{local_descent} to compute $\bbd_{i,t}^{(2)}$. Repeating this procedure $K$ times, nodes ends up having determined their NN-$K$ step component $\bbd_{i,t}^{(K)}$ .

The NN-$K$ method is summarized in Algorithm \ref{algo_NN1}. The descent iteration in \eqref{update_formula_NN} is implemented in Step 9. Implementation of this descent requires access to the NN-$K$ descent direction $\bbd_{i,t}^{(K)}$ which is computed by the loop in steps 4-8. Step $4$ initializes the loop by computing the NN-0 step $\bbd_{i,t}^{(0)}=-\bbD_{ii,t}^{-1}\bbg_{i,t}$. The core of the loop is in Step 7 which corresponds to the recursion in \eqref{local_descent}. Step 6 stands for the variable exchange that is necessary to implement Step 7. After $K$ iterations through this loop the NN-$K$ descent direction $\bbd_{i,t}^{(K)}$ is computed and can be used in Step 9. Both, steps 4 and 9, require access to the local gradient component $\bbg_{i,t}$. This is evaluated in Step 3 after receiving the prerequisite information in Step 2.

%

\subsection{Adaptive Network Newton} \label{sec:implement}

%
\begin{algorithm}[t]{
\caption{Computation of NN-$K$ step at node $i$.}
\label{algo_NNK} 
\begin{algorithmic}[1]
   \STATE \textbf{function}  
           $\bbx_{i}$
          = NN-$K$$\left(\alpha, \bbx_{i},tol \right)$ 
   \WHILE{$\|\bbg_{i}\|>tol$}
   \STATE $\bbB$ matrix blocks: $\bbB_{ii}=(1-w_{ii})\bbI$ and $\bbB_{ij}=w_{ij}\bbI$
   \STATE $\bbD$ matrix block: $\bbD_{ii,t}= \alpha \nabla^2 f_{i}(\bbx_{i}) + 2(1-w_{ii})\bbI $
   \STATE Exchange iterates $\bbx_{i}$ with neighbors $j\in \mathcal{N}_i$.
   \STATE Gradient:    
          $\displaystyle{
          \bbg_{i} = (1-w_{ii})\bbx_{i} 
                       - \sum_{j\in \mathcal{N}_i} w_{ij} \bbx_{j}
                       +\alpha \nabla f_{i}(\bbx_{i}).}$  
   \STATE Compute NN-0 descent direction $\bbd_{i}^{(0)}=-\bbD_{ii}^{-1}\bbg_{i}$\\ 
   \FOR  {$k=  0, \ldots, K-1$ } 
      \STATE Exchange elements $\bbd_{i}^{(k)}$ of the NN-$k$ step with neighbors
      \STATE NN-$(k+1)$ step:
             $\displaystyle{  
             \bbd_{i}^{(k+1)} = \bbD_{ii}^{-1}
             			\bigg[
                                   \sum_{j\in \mathcal{N}_i,j=i}\bbB_{ij} \bbd_{j}^{(k)} 
                                    - \bbg_{i}\bigg]}$.        
   \ENDFOR
          %
          \STATE Update local iterate: 
          $\displaystyle{\bbx_{i}=\bbx_{i} +\eps\ \bbd_{i}^{(K)}}$.
          \ENDWHILE
\end{algorithmic}}\end{algorithm}

As mentioned in Section \ref{sec:problem}, NN-$K$ algorithm instead of solving  \eqref{original_optimization_problem1} or its equivalent  \eqref{original_optimization_new_notation}, solves a penalty version of \eqref{original_optimization_new_notation} as introduced in \eqref{centralized_opt_problem}. 
The optimal solutions of optimization problems \eqref{original_optimization_new_notation} and \eqref{centralized_opt_problem} are different and the gap between them is upper bounded by $O(\alpha)$ \cite{YuanQing}. This observation implies that by setting a decreasing policy for $\alpha$ or equivalently an increasing policy for penalty coefficient $1/\alpha$, the solution of  \eqref{original_optimization_new_notation} approaches the minimizer of  \eqref{centralized_opt_problem}, i.e. $\tby^* \to \bby^*$ for $\alpha\to  0$.

We introduce Adaptive Network Newton-$K$ (ANN-$K$) as a version of NN-$K$ that uses a decreasing sequence of $\alpha_t$ to achieve exact convergence to the optimal solution of \eqref{original_optimization_problem1}. The idea of ANN-$K$ is to decrease parameter $\alpha_t$ by multiplying by $\eta<1$, i.e., $\alpha_{t+1}=\eta\alpha_t$, when the sequence generated by NN-$K$ is converged for a specific value of $\alpha$. To be more precise, each node $i$ has a \textit{signal vector} $\bbs_{i}=[s_{i1};\dots;s_{in}]\in \{0,1\}^n$ where each component is a binary variable. Note that $s_{ij}$ corresponds to the occurrence of receiving a signal at node $i$ from node $j$. Hence, nodes initialize their signaling components by $0$ for all the nodes in the network. At iteration $t$ node $i$ computes its local gradient norm $\|\bbg_{i,t}\|$. If the norm of gradient is smaller than a specific value called $tol$, i.e. $\|\bbg_{i,t}\|\leq tol$, it sets the local signal component to $s_{ii}=1$ and sends a signal to all the nodes in the network. The receiver nodes set the corresponding component of node $i$ in their local signal vectors to 1, i.e. $s_{ji}=1$ for $j\neq i$. This procedure implies that the signal vectors of all nodes in the network are always synchronous. The update for parameter $\alpha_t$ occurs when all the components of signal vector are 1 which is equivalent to achieving the required accuracy for all nodes in the network. Since the number of times that $\alpha_t$ should be updated is small, the cost of communication for updating $\alpha_t$ is affordable.

The ANN-$K$ method is summarized in Algorithm \ref{algo_ANN}. At each iteration of ANN-$K$ algorithm at Step 2 function NN-$K$ Step is called to update variable $\bbx_{i,t}$ for node $i$. Note that function NN-$K$ which is introduced in Algorithm \ref{algo_NNK}, runs NN-$K$ step until the time that norm of local gradient is smaller than a threshold $\|\bbg_{i}\|\leq tol$. After achieving this accuracy, in Steps 3 node $i$ updates its local signal component $s_{ii}$ to $1$ and sends it to the other nodes. In Step 4 each node $i$ updates the signal vector components of other nodes in the network. Then, in Step 6 the nodes update the penalty parameter for the next iteration as $\alpha_{t+1}=\eta\alpha_t $ if all the components of signal vector is 1, otherwise they use the previous value $\alpha_{t+1}=\alpha_{t}$. In order to reset the system after updating $\alpha_t$, all signal vectors are set to 0, i.e. $\bbs_{i}=\bb0$ for $i=1,\dots,n$ as in Step 7.

%
\begin{algorithm}[t]
\caption{Adaptive Network Newton-$K$ method at node $i$}\label{algo_ANN} 
\begin{algorithmic}[1] {
\REQUIRE  Initial iterate $\bbx_{i,0}$, initial penalty parameter $\alpha_0$ and initial sequence of bits $\bbs_i=[s_{i1};\dots;s_{in}]=[0;\dots;0].$
\FOR {$t=0,1,2,\ldots$}
     \STATE Call NN-$K$ function:
          $\bbx_{i,t+1}$
          = NN-$K$$\left(\alpha_t, \bbx_{i,t}, tol \right)$ 

   \STATE Set $s_{ii}=1$ and broadcast scalar it to all nodes.

   \STATE Set $s_{ij}=1$ for all nodes $j$ that sent a signal.

   \IF{$s_{ij}=1$ for all $j=1,\dots,n$}
   \STATE Update penalty parameter $\alpha_{t+1}=\eta\alpha_t.$
   \STATE Set $s_{ij}=0$ for all $j=1,\dots,n$.
   \ENDIF

\ENDFOR}
\end{algorithmic}\end{algorithm}


\section{Convergence Analysis}\label{sec:convergence_analysis}

In this section we show that as time progresses the sequence of objective function $F(\bby_t)$ defined in \eqref{centralized_opt_problem} approaches the optimal objective function value $F(\bby^{*})$ by considering the following assumptions.


\begin{assumption}\label{ass_weight_bounds} There exists constants $0\leq\delta<\Delta<1$ that lower and upper bound the diagonal weights for all $i$, 
\begin{equation}\label{bounds_for_local_weights}
   0\leq \delta  \leq w_{ii} \leq \Delta <1  \qquad  i=1,\ldots,n .
\end{equation}\end{assumption}


\begin{assumption}\label{convexity_assumption} 
The eigenvalues of local objective function Hessians $\nabla^2 f_{i}(\bbx)$ are bounded with positive constants $0<m\leq M<\infty$, i.e.
\begin{equation}\label{local_hessian_eigenvlaue_bounds}
m\bbI\preceq \nabla^2 f_i(\bbx)\preceq M\bbI.
\end{equation}
\end{assumption}


\begin{assumption}\label{Lipschitz_assumption} The local objective function Hessians $\nabla^2 f_i(\bbx)$ are Lipschitz continuous with parameter $L$ with respect to Euclidian norm,
\begin{equation}
   \| \nabla^2 f_i(\bbx)-\nabla^2 f_i(\hbx) \| \ \leq\  L\ \| \bbx- \hbx \|.
\end{equation}
\end{assumption}

Linear convergence of objective function $F(\bby_t)$ to the optimal objective function $F(\bby^*)$ is shown in \cite{NNICASSP} which we mention as a reference.

%
\begin{thm}\label{linear_convergence}
Consider the NN-$K$ method as defined in \eqref{diagonal_matrix}-\eqref{update_formula_NN} and the objective function $F(\bby)$ as introduced in \eqref{centralized_opt_problem}. If the stepsize $\epsilon$ is chosen as $
\epsilon = \min  \left\{ 1\ , \epsilon_0 \right\} 
$
where $\epsilon_0$ is a constant that depends on problem parameters, 
and Assumptions \ref{ass_weight_bounds}, \ref{convexity_assumption}, and \ref{Lipschitz_assumption} hold true, the sequence $F(\bby_{t})$ converges to the optimal argument $F(\bby^{*})$ at least linearly with constant $0<1-\zeta<1$. I.e.,
\begin{equation}\label{linear_convegrence_claim}
F(\bby_{t}) -F(\bby^*) \leq (1-\zeta)^t  {\left(F(\bby_{0}) -F(\bby^*) \right)}.
\end{equation}
\end{thm}


Theorem \ref{linear_convergence} shows linear convergence of sequence of objective function $F(\bby_{t})$. In the following section we study the performances of NN and ANN methods via different numerical experiments.


%
\section{Numerical analysis}\label{sec:simulations}
 
We compare the performance of DGD and different versions of NN in the minimization of a distributed quadratic objective. The comparison is done in terms of both, number of iterations and number of information exchanges. Specifically, for each agent $i$ we consider a positive definite diagonal matrix $\bbA_i\in \mbS_{p}^{++}$ and a vector $\bbb_i \in \mbR^{p}$ to define the local objective function $f_i(\bbx)  :=  ({1}/{2}) \bbx^{T} \bbA_i\bbx +\bbb_i^T\bbx$. Therefore, the global cost function $f(\bbx)$ is written as
\begin{equation}\label{example_problem}
  f(\bbx)  := \sum_{i=1}^n  \frac{1}{2} \bbx^{T} \bbA_i \bbx +\bbb_i^{T}\bbx   \ .
\end{equation}
The difficulty of solving \eqref{example_problem} is given by the condition number of the matrices $\bbA_i$. To adjust condition numbers we generate diagonal matrices $\bbA_i$ with random diagonal elements $a_{ii}$. The first $p/2$ diagonal elements $a_{ii}$ are drawn uniformly at random from the discrete set $\{1, 10^{-1},\ldots, 10^{-\xi}\}$ and the next $p/2$ are uniformly and randomly chosen from the set  $\{1,10^1,\ldots, 10^{\xi}\}$. This choice of coefficients yields local matrices $\bbA_i$ with eigenvalues in the interval $[10^{-\xi}, 10^{\xi}]$ and global matrices $\sum_{i=1}^{n}\bbA_i$ with eigenvalues in the interval $[n10^{-\xi}, n10^{\xi}]$. The condition numbers are typically $10^{2\xi}$ for the local functions and $10^\xi$ for the global objectives. The linear terms $\bbb_i^{T}\bbx$ are added so that the different local functions have different minima. The vectors $\bbb_i$ are chosen uniformly at random from the box $[0,1]^p$. 

For the quadratic objective in \eqref{example_problem} we can compute the optimal argument $\bbx^*$ in closed form. We then evaluate convergence through the relative error that we define as the average normalized squared distance between local vectors $\bbx_i$ and the optimal decision vector $\bbx^*$,
\begin{equation}\label{distance_error}
e_t := \frac{1}{n} \sum_{i=1}^{n}\frac{\| \bbx_{i,t}-\bbx^*\|^2}{\|\bbx^*\|^2}.
\end{equation}
The network connecting the nodes is a $d$-regular cycle where each node is connected to exactly $d$ neighbors and $d$ is assumed even. The graph is generated by creating a cycle and then connecting each node with the $d/2$ nodes that are closest in each direction. The diagonal weights in the matrix $\bbW$ are set to $w_{ii} = 1/2+1/2(d+1)$ and the off diagonal weights to $w_{ij}= 1/2(d+1)$ when $j\in\ccalN_i$. 

%
\begin{figure}[t]
\centering
\includegraphics[width=\linewidth,height=0.63\linewidth]{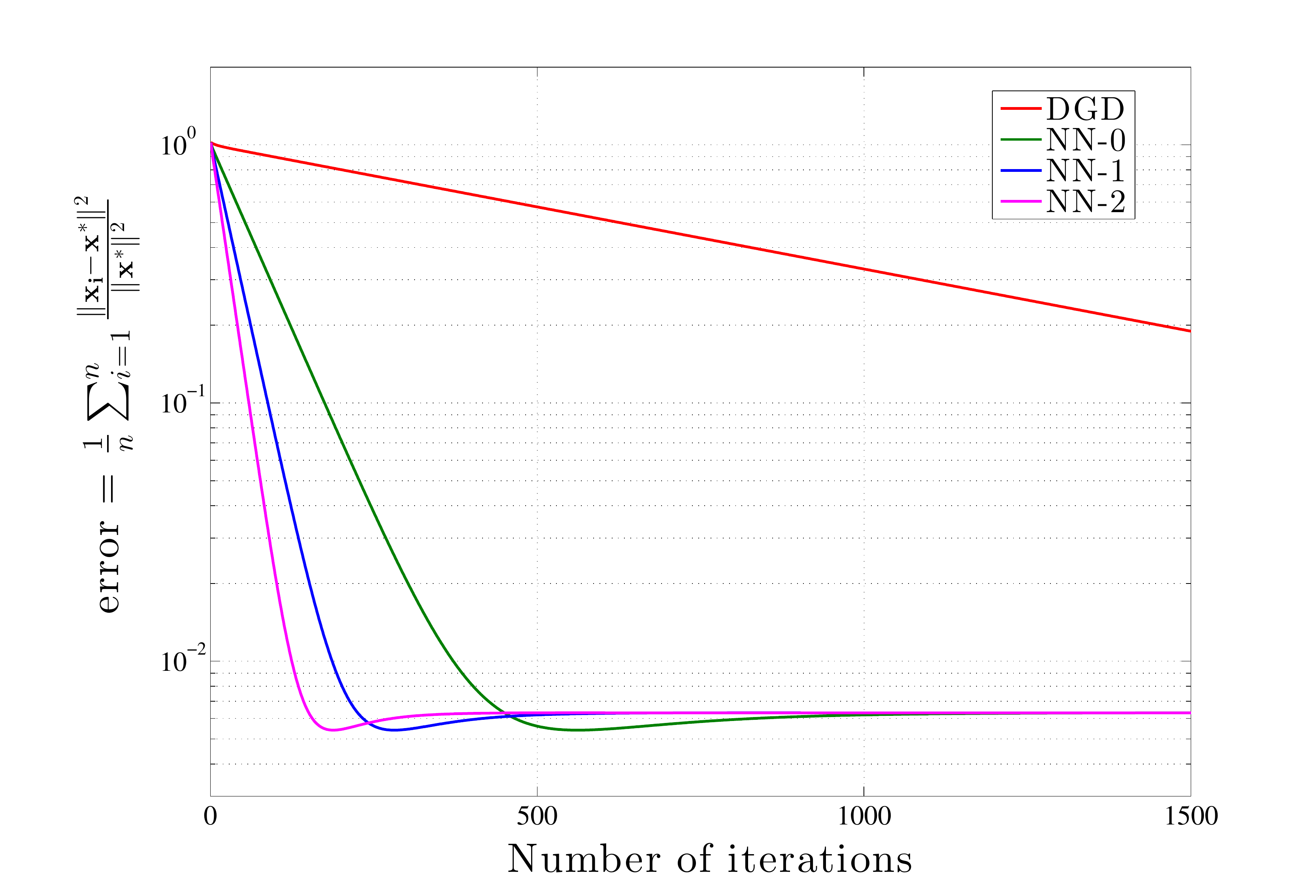}
\caption{Convergence of DGD, NN-0, NN-1, and NN-2 in terms of number of iterations. The NN methods converges faster than DGD. Furthermore, the larger $K$ is, the faster NN-$K$ converges.}
\label{fig:iter_illus}
\end{figure}

%
\begin{figure}[t]
\centering
\includegraphics[width=\linewidth,height=0.63\linewidth]{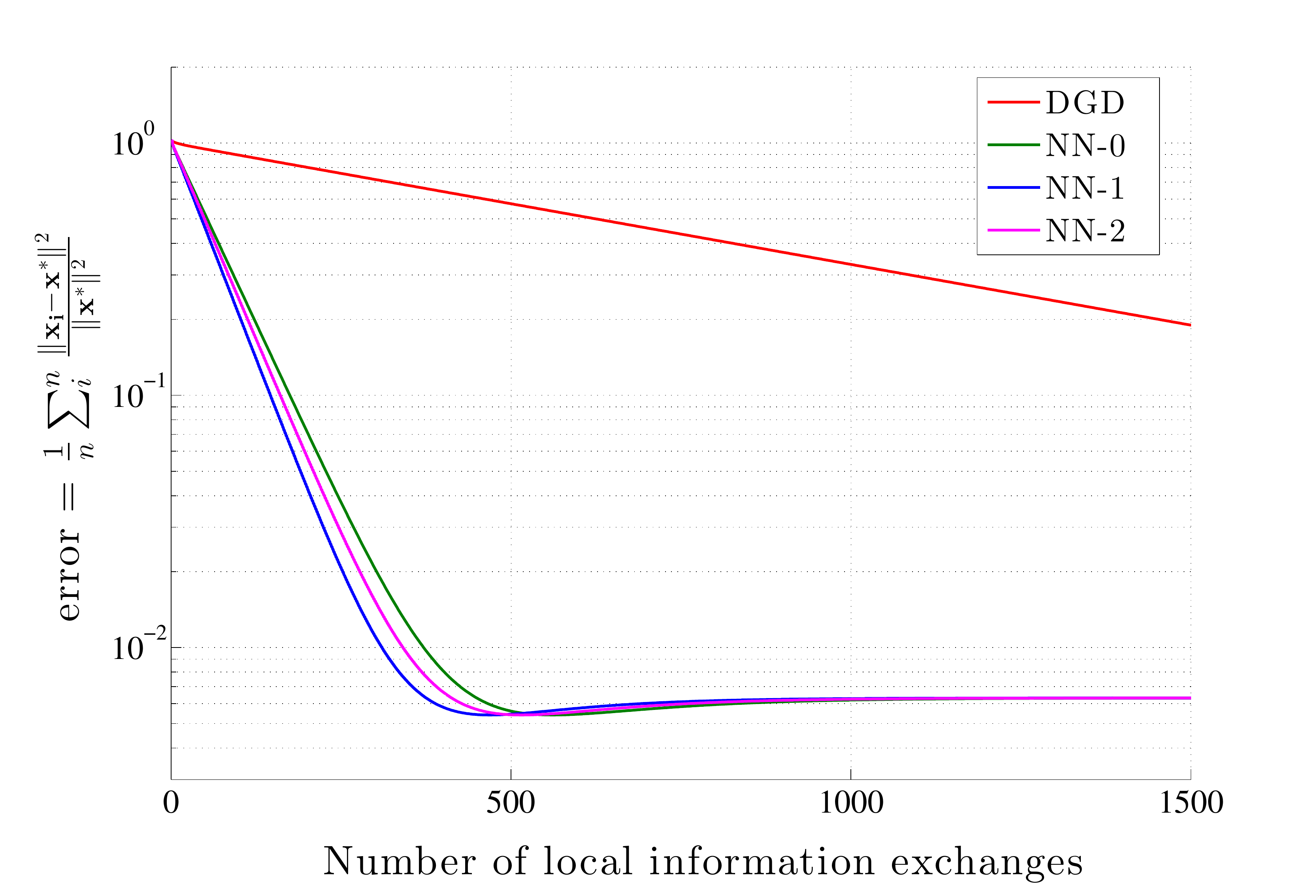}
\caption{Convergence of DGD, NN-0, NN-1, and NN-2 in terms of number of communication exchanges. The NN-$K$ methods retain the advantage over DGD but increasing $K$ may not result in faster convergence. For this particular instance it is actually NN-1 that converges fastest in terms of number of communication exchanges. }
\label{fig:local_exch_illus}
\end{figure}

%
\begin{figure}[t]
\centering
\includegraphics[width=\linewidth,height=0.63\linewidth]{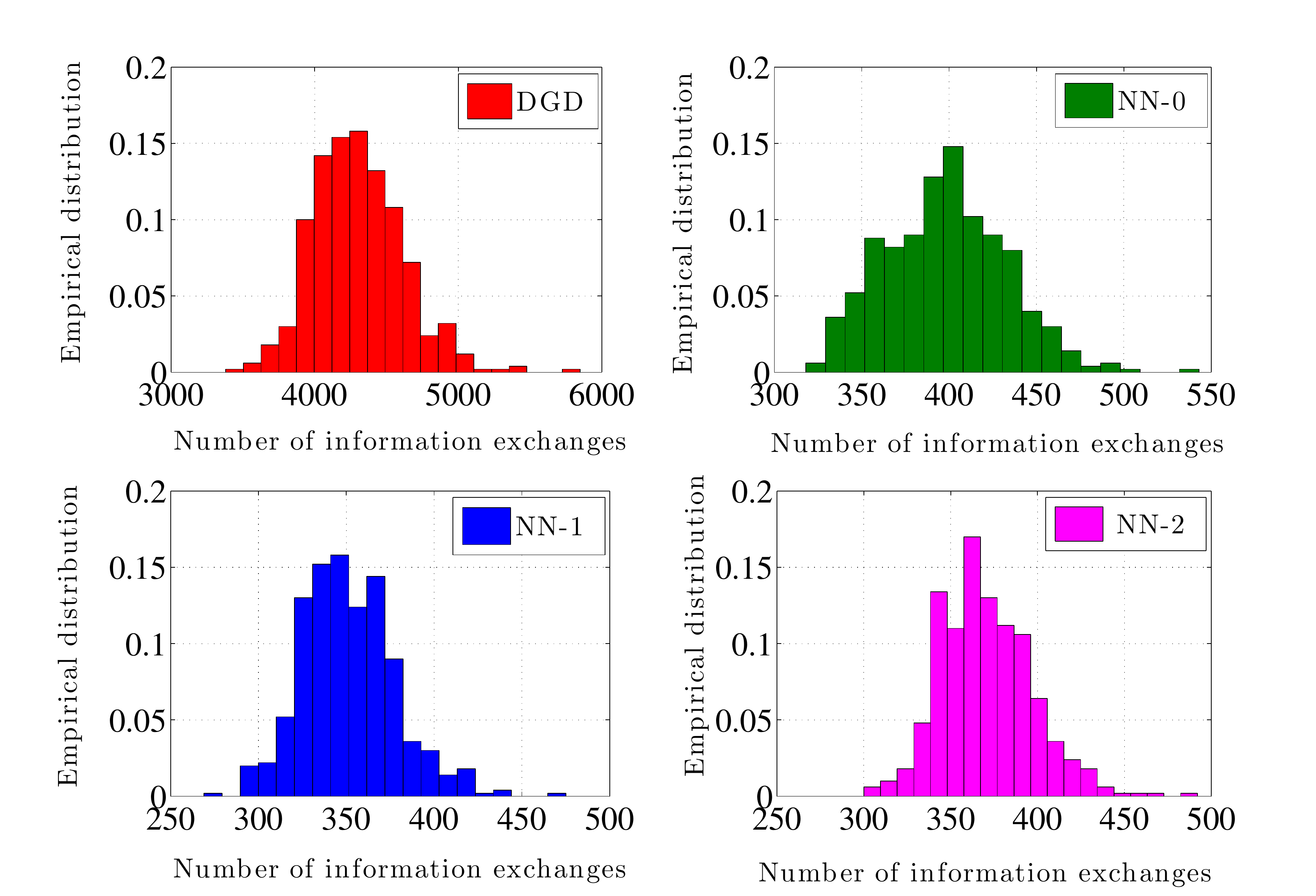}\vspace{-2mm}
\caption{Histograms of the number of information exchanges required to achieving accuracy $e_t<10^{-2}$. The qualitative observations made in figures \ref{fig:iter_illus} and \ref{fig:local_exch_illus} hold over a range of random problem realizations.}
\label{fig:local_exch_dist}
\end{figure}

%
In the subsequent experiments we set the network size to $n=100$, the dimension of the decision vectors to $p=4$, the condition number parameter to $\xi=2$, the penalty coefficient inverse to $\alpha =10^{-2}$, and the network degree to $d=4$. The NN step size is set to $\eps=1$, which is always possible when we have quadratic objectives. Figure \ref{fig:iter_illus} illustrates a sample convergence path for DGD, NN-0, NN-1, and NN-2 by measuring the relative error $e_t$ in \eqref{distance_error} with respect to the number of iterations $t$. As expected for a problem that doesn't have a small condition number -- in this particular instantiation of the function in \eqref{example_problem} the condition number is $95.2$ -- different versions of NN are much faster than DGD. E.g., after $t=1.5\times 10^3$ iterations the error associated which DGD iterates is $e_t\approx1.9\times 10^{-1}$. Comparable or better accuracy $e_t<1.9\times 10^{-1}$ is achieved in $t=132$, $t=63$, and $t=43$ iterations for NN-0, NN-1, and NN-2, respectively. 

Further recall that $\alpha$ controls the difference between the actual optimal argument  $\tby^*=[\bbx^*;\ldots;\bbx^*]$ [cf. \eqref{original_optimization_new_notation}] and the argument $\bby^*$ [cf. \eqref{centralized_opt_problem}] to which DGD and NN converge. Since we have $\alpha =10^{-2}$ and the difference between these two vectors is of order $O(\alpha)$, we expect the error in \eqref{distance_error} to settle at $e_t\approx10^{-2}$. The error actually settles at $e_t\approx6.3\times 10^{-3}$ and it takes all three versions of NN less than $t=400$ iterations to do so. It takes DGD more than $t=10^4$ iterations to reach this value. This relative performance difference decreases if the problem has better conditioning but can be made arbitrarily large by increasing the condition number of the matrix $\sum_{i=1}^{n}\bbA_i$. The number of iterations required for convergence can be further decreased by considering higher order approximations in \eqref{Hessian_approximation_iteration}. The advantages would be misleading because they come at the cost of increasing the number of communications required to approximate the Newton step.

%
To study this latter effect we consider the relative performance of DGD and different versions of NN in terms of the number of local information exchanges. Note that each iteration in NN-$K$ requires a total of $K+1$ information exchanges with each neighbor, as opposed to the single variable exchange required by DGD. After $t$ iterations the number of variable exchanges between each pair of neighbors is $t$ for DGD and $(K+1)t$ for NN-$K$. Thus, we can translate Figure \ref{fig:iter_illus} into a path in terms of number of communications by scaling the time axis by $(K+1)$. The result of this scaling is shown in Figure \ref{fig:local_exch_illus}. The different versions of NN retain a significant, albeit smaller, advantage with respect to DGD. Error $e_t<10^{-2}$ is achieved by NN-0, NN-1, and NN-2 after $(K+1)t=3.7\times 10^2$, $(K+1)t=3.1\times 10^2$, and $(K+1)t=3.4\times 10^2$ variable exchanges, respectively. When measured in this metric it is no longer true that increasing $K$ results in faster convergence. For this particular problem instance it is actually NN-1 that converges fastest in terms of number of communication exchanges.

%
For a more more comprehensive evaluation we consider $10^3$ different random realizations of \eqref{example_problem} where we also randomize the degree $d$ of the $d$-regular graph that we choose from the even numbers in the set $[2,10]$. The remaining parameters are the same used to generate figures \ref{fig:iter_illus} and \ref{fig:local_exch_illus}. For each joint random realization of network and objective we run DGD, NN-0, NN-1, and NN-2, until achieving error $e_t<10^{-2}$ and record the number of communication exchanges that have elapsed -- which amount to simply $t$ for DGD and $(K+1)t$ for NN. The resulting histograms are shown in Figure \ref{fig:local_exch_dist}. The mean times required to reduce the error to $e_t<10^{-2}$ are $4.3\times 10^3$ for DGD and $4.0\times 10^2$, $3.5\times 10^2$, and $3.7\times 10^2$ for NN-0, NN-1, and NN-2. As in the particular case shown in figures \ref{fig:iter_illus} and \ref{fig:local_exch_illus}, NN-1 performs best in terms of communication exchanges. Observe, however, that the number of communication exchanges required by NN-2 is not much larger and that NN-2 requires less computational effort than NN-1 because the number of iterations $t$ is smaller.

%

%
\begin{figure}[t]
\centering
\includegraphics[width=\linewidth,height=0.63\linewidth]{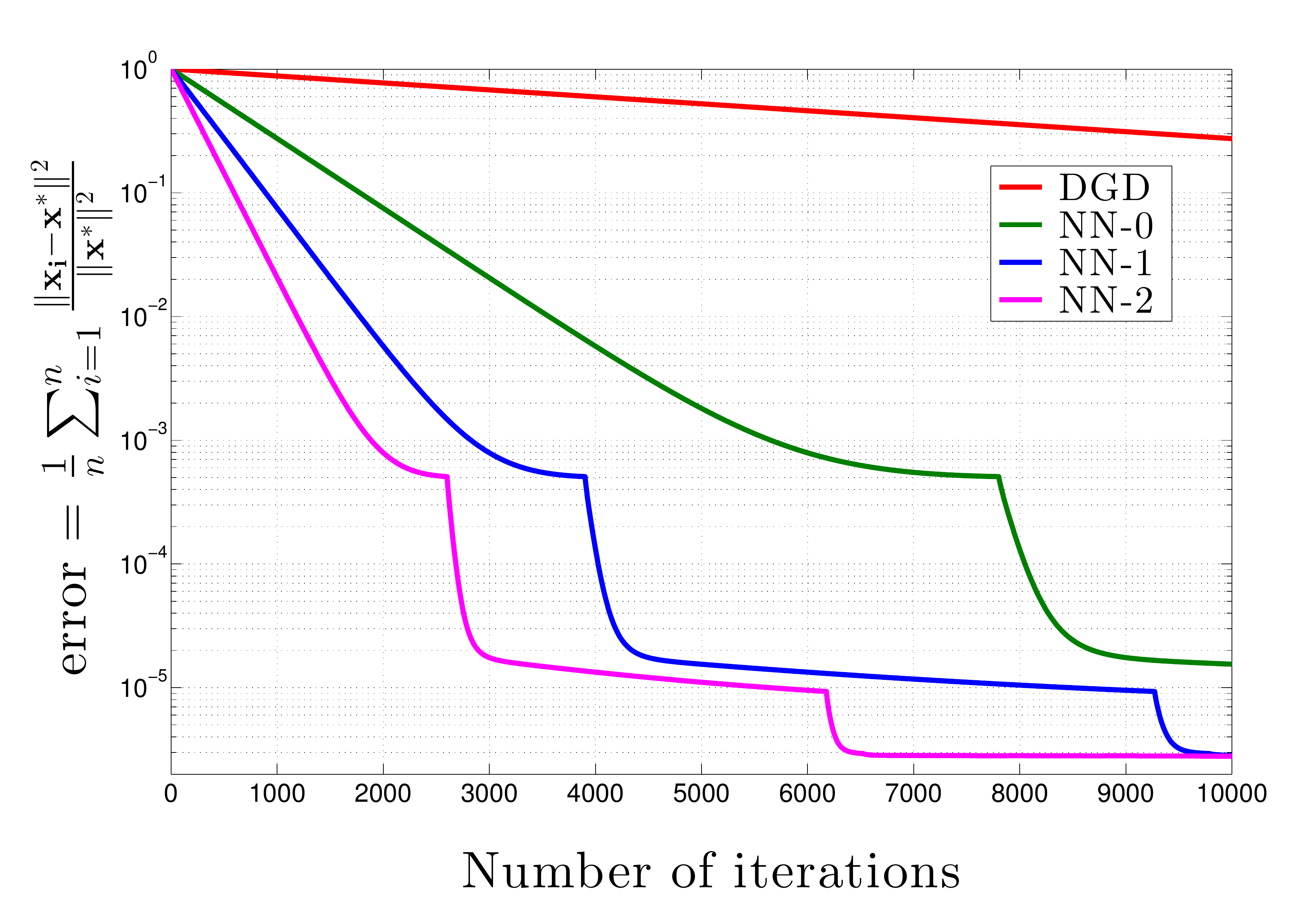}
\caption{{Convergence of adaptive DGD, NN-0, NN-1, and NN-2 for $\alpha_0\!=\!10^{-2}$.} }
\label{fig:4_2}
\end{figure}

\subsection{Adaptive Network Newton}

Given that DGD and NN are penalty methods it is of interest to consider their behavior when the inverse penalty parameter $\alpha$ is decreased recursively. The adaptation of $\alpha$ for NN-$K$ is discussed in Section \ref{sec:implement} where it is termed adaptive (A)NN-$K$. The same adaptation strategy is considered here for DGD. The parameter $\alpha$ is kept constant until the local gradient components $\bbg_{i,t}$ become smaller than a given tolerance $\text{tol}$, i.e., until $\|\bbg_{i,t}\|\leq {tol}$ for all $i$. When this tolerance is achieved, the parameter $\alpha$ is scaled by a factor $\eta<1$, i.e., $\alpha$ is decreased from its current value to $\eta\alpha$. This requires the use of a signaling method like the one summarized in Algorithm \ref{algo_ANN} for ANN-$K$.

\begin{figure}[t]
\centering
\includegraphics[width=\linewidth,height=0.63\linewidth]{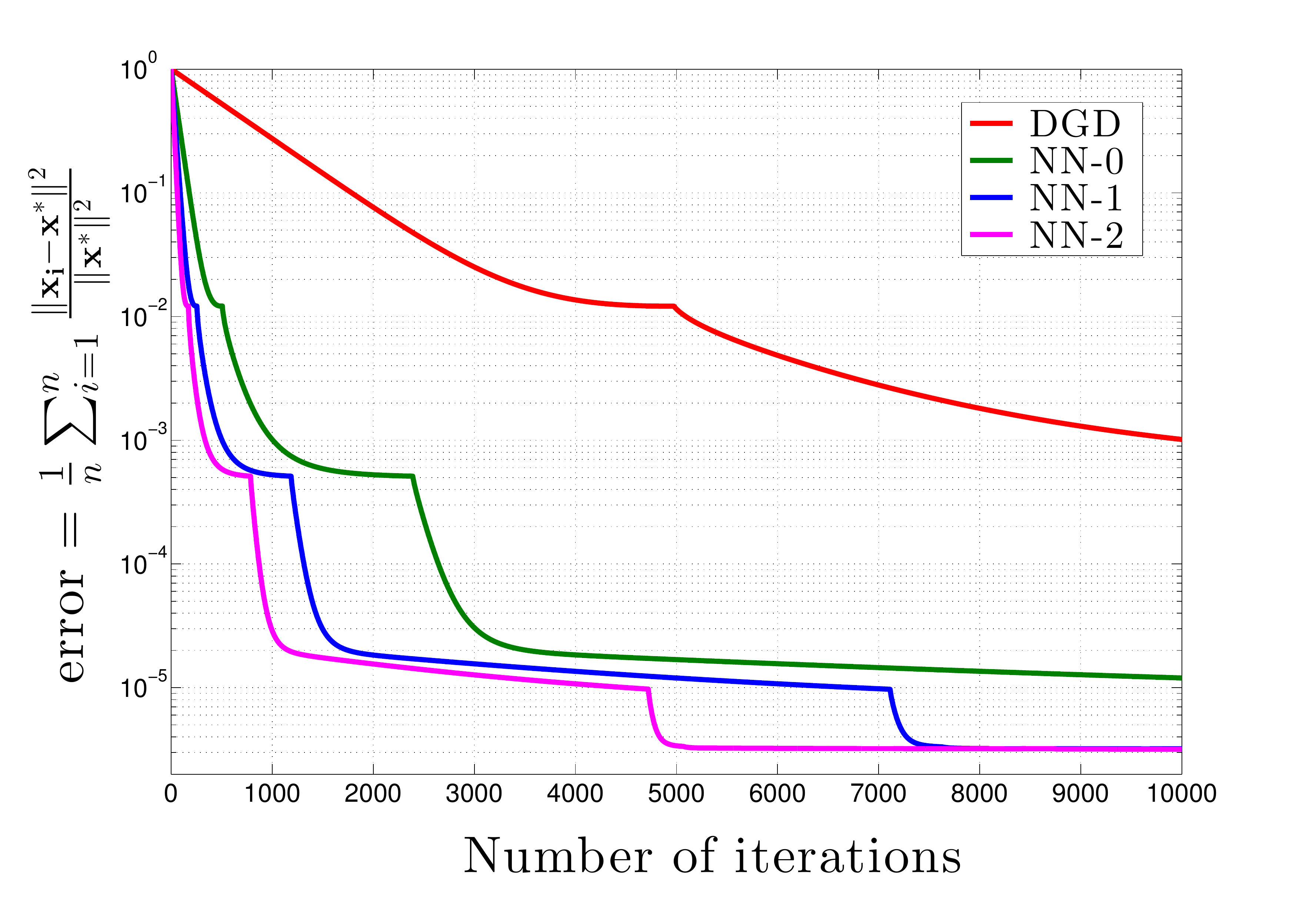}
\caption{{\!Convergence of Adaptive DGD, NN-0, NN-1, and NN-2 for $\alpha_0\!=\!10^{-1}$.} }
\label{fig:4}
\end{figure}

We consider the objective in \eqref{example_problem} and nodes connected by a $d$-regular cycle. We use the same parameters used to generate figures \ref{fig:iter_illus} and \ref{fig:local_exch_illus}. The adaptive gradient tolerance is set to ${tol}=10^{-3}$ and the scaling parameter to $\eta=0.1$. We consider two different scenarios where the initial penalty parameters are $\alpha=\alpha_{0}=10^{-1}$ and $\alpha=\alpha_0=10^{-2}$. The respective error trajectories $e_t$ with respect to the number o iterations are shown in figures \ref{fig:4_2} -- where  $\alpha_0=10^{-2}$ -- and \ref{fig:4} -- where $\alpha_0=10^{-1}$. In each figure we show $e_t$ for adaptive DGD, ANN-0, ANN-1, and ANN-2. Both figures show that the ANN methods outperform adaptive DGD and that larger $K$ reduces the number of iterations that it takes ANN-$K$ to achieve a target error. These results are consistent with the findings summarized in figures \ref{fig:iter_illus}-\ref{fig:local_exch_dist}. 

More interesting conclusions follow from a comparison across figures \ref{fig:iter_illus} and \ref{fig:local_exch_illus}. We can see that it is better to start with the (larger) value $\alpha=10^{-1}$ even if the method initially converges to a point farther from the actually optimum. This happens because problems with larger $\alpha$ are better conditioned and thus easier to minimize.

%

\bibliographystyle{IEEEtran}
  \bibliography{bmc_article}
   \end{document}